\documentclass{amsart}
\usepackage{amssymb,latexsym}
\usepackage{amscd,amsthm}
\usepackage{tikz-cd}

\usepackage[all]{xy}

\newtheorem{theorem}{Theorem}[section]

\theoremstyle{definition}

\DeclareMathOperator{\Ext}{Ext}

\newcommand{\cat}[1]{\mathcal{#1}}           

\newcommand{\class}[1]{\mathcal{#1}}   

\newcommand{\tilclass}[1]{\widetilde{\class{#1}}}

\begin{document}

\title{How to construct a Hovey triple from two cotorsion pairs}

\author{James Gillespie}
\address{Ramapo College of New Jersey \\
         School of Theoretical and Applied Science \\
         505 Ramapo Valley Road \\
         Mahwah, NJ 07430}
\email[Jim Gillespie]{jgillesp@ramapo.edu}
\urladdr{http://pages.ramapo.edu/~jgillesp/}

\date{\today}

\begin{abstract}
Let $\cat{A}$ be an abelian category, or more generally a weakly idempotent complete exact category, and suppose we have two complete hereditary cotorsion pairs $(\class{Q}, \tilclass{R})$ and $(\tilclass{Q}, \class{R})$ in $\cat{A}$ satisfying $\tilclass{R} \subseteq \class{R}$ and $\class{Q} \cap \tilclass{R} = \tilclass{Q} \cap \class{R}$. We show how to construct a (necessarily unique) abelian model structure on $\cat{A}$ with $\class{Q}$ (respectively $\tilclass{Q}$) as the class of cofibrant (resp. trivially cofibrant) objects and $\class{R}$ (respectively $\tilclass{R}$) as the class of fibrant (resp. trivially fibrant) objects.
\end{abstract}

\maketitle

\section{Introduction}\label{sec-intro}

Let $\cat{A}$ be an abelian category. In~\cite{hovey} we learned a one-to-one correspondence between complete cotorsion pairs in $\cat{A}$ and abelian model structures on $\cat{A}$. We assume the reader is familiar with this correspondence. Let's recall that an abelian model structure on $\cat{A}$ is equivalent to a triple of classes $(\class{Q},\class{W},\class{R})$ for which $\class{W}$ is thick and which admits two complete cotorsion pairs $(\class{Q},\class{W} \cap \class{R})$ and $(\class{Q} \cap \class{W},\class{R})$. This data is equivalent to an abelian model structure on $\cat{A}$ where $\class{W}$ are the trivial objects, $\class{Q}$ are the cofibrant objects and $\class{R}$ are the fibrant objects. In practice virtually all the abelian model structures we encounter also have the extra condition that they are \emph{hereditary}. This just means that $\class{Q}$ is closed under taking kernels of epimorphisms between objects in $\class{Q}$ and $\class{R}$ is closed between taking cokernels of monomorphisms between objects in $\class{R}$. We will assume this throughout this note, that the cotorsion pairs are hereditary.  Now let's denote these cotorsion pairs by $$(\class{Q}, \tilclass{R}) = (\class{Q},\class{W} \cap \class{R}) \text{ and } (\tilclass{Q}, \class{R}) = (\class{Q} \cap \class{W},\class{R}).$$
It is easy to note the following properties of these two complete hereditary cotorsion pairs:
\begin{enumerate}
\item $\tilclass{R} \subseteq \class{R}$ and $\tilclass{Q} \subseteq \class{Q}$.
\item $\tilclass{Q} \cap \class{R} = \class{Q} \cap \tilclass{R}$.
\end{enumerate}
Remarkably, there is a converse to this as we will prove in this note. So our goal is to prove the following theorem.

\begin{theorem}\label{them-main}
Let $(\class{Q}, \tilclass{R})$ and $(\tilclass{Q}, \class{R})$ be two complete hereditary cotorsion pairs satisfying conditions \textnormal{(1)} and \textnormal{(2)} above. Then there is a unique thick class $\class{W}$ for which $(\class{Q},\class{W},\class{R})$ is a Hovey triple. Moreover, this thick class $\class{W}$ can be described in the two following ways:
\begin{align*}
   \class{W}  &= \{\, X \in \class{A} \, | \, \exists \, \text{a short exact sequence } \, X \rightarrowtail R \twoheadrightarrow Q \, \text{ with} \, R \in \tilclass{R} \, , Q \in \tilclass{Q} \,\} \\
           &= \{\, X \in \class{A} \, | \, \exists \, \text{a short exact sequence } \, R' \rightarrowtail Q' \twoheadrightarrow X \, \text{ with} \, R' \in \tilclass{R} \, , Q' \in \tilclass{Q} \,\}.
          \end{align*}
\end{theorem}

The proof we give is both elementary and very general. It holds in the general setting of when $\cat{A}$ is a weakly idempotent complete exact category. The author showed in~\cite{gillespie-exact model structures} that Hovey's correspondence carries over to this setting. So while we state the theorems in the setting of abelian categories we have particularly constructed the proofs so that they hold in a general weakly idempotent complete exact category.

The author again wishes to thank Hanno Becker and the referee of his paper~\cite{gillespie-recollements}. The construction of Hovey triples that we give here is a direct generalization of a construction of Becker from~\cite{becker}. In particular, when the two given cotorsion pairs are \emph{injective}, then our construction is exactly Becker's right Bousfield localization construction from~\cite{becker}.  On the other hand, when the two given cotorsion pairs are \emph{projective}, then our construction coincides with Becker's left Bousfield localization construction from~\cite{becker}. The main difference in our proof, when comparing it to Becker's proof, is in how we show that the class of trivial objects $\class{W}$ is thick. Our proof, while longer, is more direct and uses only elementary properties of short exact sequences.

Finally, the author wishes to thank Mark Hovey, who the author has been fortunate to work with over the years. It was during a recent exchange that the author found Theorem~\ref{them-main}, and it solves the problem of finding the Gorenstein AC-flat model structure on the category of (left) modules over an arbitrary ring. This is the flat analog of the Gorenstein AC-injective and Gorenstein AC-projective model structures of~\cite{bravo-gillespie-hovey}. These applications, and others, will appear elsewhere. But the author finds Theorem~\ref{them-main} to be of interest in its own right.

\section{Proof of the theorem}\label{sec-proof}

We wish to prove Theorem~\ref{them-main}. So assume $(\class{Q}, \tilclass{R})$ and $(\tilclass{Q}, \class{R})$ are two complete hereditary cotorsion pairs satisfying conditions
\begin{enumerate}
\item $\tilclass{R} \subseteq \class{R}$ and $\tilclass{Q} \subseteq \class{Q}$.
\item $\tilclass{Q} \cap \class{R} = \class{Q} \cap \tilclass{R}$.
\end{enumerate}
We wish to construct a Hovey triple $(\class{Q},\class{W},\class{R})$ with the properties in Theorem~\ref{them-main}.

\begin{proof}
We start by showing that the two classes below that define $\class{W}$ do coincide.
\begin{align*}
   \class{W}  &= \{\, X \in \class{A} \, | \, \exists \, \text{a short exact sequence } \, X \rightarrowtail R \twoheadrightarrow Q \, \text{ with} \, R \in \tilclass{R} \, , Q \in \tilclass{Q} \,\} \\
           &= \{\, X \in \class{A} \, | \, \exists \, \text{a short exact sequence } \, R' \rightarrowtail Q' \twoheadrightarrow X \, \text{ with} \, R' \in \tilclass{R} \, , Q' \in \tilclass{Q} \,\}.
          \end{align*}
So say $X$ is in the top class, that is, that there is a short exact sequence $X \rightarrowtail R \twoheadrightarrow Q$ where $R \in \tilclass{R}$ and $Q \in \tilclass{Q}$. Since $(\class{Q},\tilclass{R})$ has enough projectives we can find a short exact sequence $R' \rightarrowtail Q' \twoheadrightarrow R$ where $R' \in \tilclass{R}$ and $Q' \in \class{Q}$. We take a pullback and from~\cite[Proposition~2.12]{buhler-exact categories} we get the commutative diagram below with exact rows and columns and whose lower left corner is a bicartesian (pushpull) square.
$$\begin{tikzcd}
R' \arrow[equal]{r} \arrow[tail]{d} & R' \arrow[tail]{d} &  \\
P \arrow[tail]{r} \arrow[two heads]{d} & Q' \arrow[two heads]{r} \arrow[two heads]{d} & Q \arrow[equal]{d} \\
X \arrow[tail]{r} & R \arrow[two heads]{r} & Q
\end{tikzcd}$$
Since $\tilclass{R}$ is closed under extensions, by the assumption (2) $\tilclass{Q} \cap \class{R} = \class{Q} \cap \tilclass{R}$, we deduce that $Q' \in  \tilclass{Q} \cap \class{R}$. Now since $Q',Q \in \tilclass{Q}$ and the cotorsion pairs are hereditary we conclude $P \in \tilclass{Q}$. Now the left vertical column shows that $X$ is in the bottom class describing $\class{W}$. A similar argument will show that any $X$ in the bottom class must be in the top class. So the two descriptions of $\class{W}$ coincide.

\

\noindent ($\class{W}$ is Thick.)  We must show $\class{W}$ is closed under retracts and that whenever two out of three terms in a short exact sequence $X  \rightarrowtail Y \twoheadrightarrow  Z$ are in $\class{W}$ then so is the third.

We start by showing $\class{W}$ is closed under retracts. So let $W \in \class{W}$ and $X \xrightarrow{i} W \xrightarrow{p} X$ be such that $pi = 1_X$. We wish to show $X \in \class{W}$. Start by writing $W  \rightarrowtail \widetilde{R} \twoheadrightarrow  \widetilde{Q}$ and also apply the fact that $(\tilclass{Q}, \class{R})$ has enough injectives to get a short exact sequence $X  \rightarrowtail R \twoheadrightarrow  \widetilde{Q}'$. Now construct a commutative diagram as shown below:
$$\begin{CD}
X       @>>>     R      @>>>      \widetilde{Q}' \\
@V i VV             @VV j V           @VVV     \\
W       @>>> \widetilde{R} @>>>    \widetilde{Q} \\
@V p VV            @VV q V           @VVV\\
X       @>>>       R      @>>>     \widetilde{Q}' \\
\end{CD}$$
The map $j$ exists because $\Ext^1_{\class{A}}(\widetilde{Q}', \widetilde{R}) = 0$ and similarly the map $q$ exists because $\Ext^1_{\class{A}}(\widetilde{Q}, R) = 0$. Next, the two right vertical maps exist simply by the universal property of $\widetilde{Q}'$ and $\widetilde{Q}$ being cokernels. Next denote the map $X  \rightarrowtail R$ by $k$, and its cokernel $R \twoheadrightarrow  \widetilde{Q}'$ by $h$. Then we see $(1_{R} - qj)k = k - qjk = k - k = 0$. So again the universal property of $\widetilde{Q}'$ being the cokernel of $k$ gives us a map $\widetilde{Q}' \xrightarrow{t} R$ such that $th = 1_{R} - qj$. This proves that $1_{R} - qj$ factors through $\widetilde{Q}'$,  but now we argue that this implies $1_{R} - qj$ actually factors through an object of $\tilclass{Q} \cap \class{R} = \class{Q} \cap \tilclass{R}$. Indeed using that $(\tilclass{Q}, \class{R})$ has enough injectives we find $\widetilde{Q}'  \rightarrowtail R' \twoheadrightarrow  \widetilde{Q}''$, but this time it follows that $R' \in \tilclass{Q} \cap \class{R}$.  Now using that $\Ext^1_{\class{A}}(\widetilde{Q}'', R) = 0$, we see that $\widetilde{Q}' \xrightarrow{t} R$ extends over $\widetilde{Q}'  \rightarrowtail R'$. So we see that we can find maps $R \xrightarrow{\alpha} R'$ and $R' \xrightarrow{\beta} R$ such that $1_{R} - qj = \beta \alpha$ and $R' \in \tilclass{Q} \cap \class{R} = \class{Q} \cap \tilclass{R}$. Thus the composition below is the identity $1_{R}$.
$$R \xrightarrow{(j \ \alpha)} \widetilde{R} \oplus R' \xrightarrow{q+\beta} R$$ But this just means $R$ is a retract of $\widetilde{R} \oplus R'$. Since $\tilclass{R}$ is closed under direct sums and retracts get that $R \in \tilclass{R}$. This proves $X \in \class{W}$ and we are done.

\

We now turn to the two out of three property and our next immediate goal is to show closure of $\class{W}$ under extensions.  Note now that $\class{W}$ clearly contains both $\tilclass{Q}$ and $\tilclass{R}$, for we will use this ahead. We start by making the following claim.

\

\noindent \underline{Claim 1}. Suppose $R  \rightarrowtail Y \twoheadrightarrow  W$ is exact with $R \in \tilclass{R}$ and $W \in \class{W}$. Then there exists a commutative diagram as below where $\widetilde{Q}, \widetilde{Q}', \widetilde{Q}'' \in \tilclass{Q}$ and $\widetilde{R}, \widetilde{R}', \widetilde{R}'' \in \tilclass{R}$.
$$\begin{tikzcd}
\widetilde{R} \arrow[tail]{r} \arrow[tail]{d} & \widetilde{R}'' \arrow[two heads]{r} \arrow[tail]{d} & \widetilde{R}' \arrow[tail]{d} \\
\widetilde{Q} \arrow[tail]{r} \arrow[two heads]{d} & \widetilde{Q}'' \arrow[two heads]{r} \arrow[two heads]{d} & \widetilde{Q}' \arrow[two heads]{d} \\
R \arrow[tail]{r}  & Y \arrow[two heads]{r} & W
\end{tikzcd}$$ Indeed since $R, W \in \class{W}$ there are short exact sequence $\widetilde{R}  \rightarrowtail \widetilde{Q} \twoheadrightarrow  R$ and $\widetilde{R}'  \rightarrowtail \widetilde{Q}' \twoheadrightarrow  W$ with $\widetilde{R} , \widetilde{R}' \in \tilclass{R}$ and $\widetilde{Q}, \widetilde{Q}' \in \tilclass{Q}$. But in this case we also have $R \in \tilclass{R} \subseteq \class{R}$ and so $\Ext^1_{\cat{A}}(\widetilde{Q}',R) = 0$. This means there exists a lift as shown:
$$\begin{tikzcd}
&  & \widetilde{Q}' \arrow[two heads]{d} \arrow[dashed]{dl} \\
R \arrow[tail]{r}  & Y \arrow[two heads]{r} & W
\end{tikzcd}$$
This lift allows for the construction, analogous to the usual Horseshoe Lemma of homological algebra, of a commutative diagram as below.
$$\begin{tikzcd}
\widetilde{R} \arrow[tail]{r} \arrow[tail]{d} & \widetilde{R}'' \arrow[two heads]{r} \arrow[tail]{d} & \widetilde{R}' \arrow[tail]{d} \\
\widetilde{Q} \arrow[tail]{r} \arrow[two heads]{d} & \widetilde{Q} \oplus \widetilde{Q}' \arrow[two heads]{r} \arrow[two heads]{d} & \widetilde{Q}' \arrow[two heads]{d} \\
R \arrow[tail]{r}  & Y \arrow[two heads]{r} & W
\end{tikzcd}$$
Since any class that is part of a cotorsion pair is closed under direct sums and extensions we now have  $\widetilde{Q} \oplus \widetilde{Q}' \in \tilclass{Q}$ and $\widetilde{R}'' \in \tilclass{R}$, and so we have proved our first claim.

\

\noindent \underline{Claim 2}. $\class{W}$ is closed under extensions. Now suppose $W  \rightarrowtail Y \twoheadrightarrow  W'$ is exact with $W, W' \in \class{W}$. We need to prove that $Y \in \class{W}$ too. Since $W \in \class{W}$ we may now find an exact sequence $W  \rightarrowtail R \twoheadrightarrow  Q$ where $R \in \tilclass{R}$ and $Q \in \tilclass{Q}$. Now form the pushout diagram shown.
$$\begin{tikzcd}
W \arrow[tail]{r} \arrow[tail]{d} & Y \arrow[two heads]{r} \arrow[tail]{d} & W' \arrow[equal]{d} \\
R \arrow[tail]{r} \arrow[two heads]{d} & P \arrow[two heads]{r} \arrow[two heads]{d} & W'  \\
Q \arrow[equal]{r} & Q
\end{tikzcd}$$ Note the second row is the type of row from claim 1. So the Horseshoe argument provides a diagram as shown below where $\widetilde{Q}, \widetilde{Q}', \widetilde{Q}'' \in \tilclass{Q}$ and $\widetilde{R}, \widetilde{R}', \widetilde{R}'' \in \tilclass{R}$.
$$\begin{tikzcd}[row sep=scriptsize, column sep=scriptsize]
& & & & & \widetilde{R}' \arrow[tail]{dd}  \\
& & & \widetilde{R}'' \arrow[tail]{dd} \arrow[two heads]{urr} & &  \\
& \widetilde{R} \arrow[tail]{dd} \arrow[tail]{urr} & & & & \widetilde{Q}' \arrow[two heads]{dd} \\
& & & \widetilde{Q}'' \arrow[two heads]{dd}  \arrow[two heads]{urr} & W' \arrow[equal]{dr} & \\
 & \widetilde{Q} \arrow[two heads]{dd} \arrow[tail]{urr} & Y  \arrow[tail]{dr} \arrow[two heads]{urr} & & & W'  \\
W \arrow[tail]{dr}  \arrow[tail]{urr} & & & P \arrow[two heads]{urr}& & \\
& R  \arrow[tail]{urr}& & & & \\
\end{tikzcd}$$
Now pullback the entire diagram over the original exact sequence $W  \rightarrowtail Y \twoheadrightarrow  W'$ and we get what we want. In particular, the pullback in the middle of the diagram leads to the following bicartesian (pushpull) square and the hereditary property of $\tilclass{Q}$ gives us $L \in \tilclass{Q}$.
$$\begin{tikzcd}
\widetilde{R}'' \arrow[equal]{r} \arrow[tail]{d} & \widetilde{R}'' \arrow[tail]{d} &  \\
L \arrow[tail]{r} \arrow[two heads]{d} & \widetilde{Q}'' \arrow[two heads]{r} \arrow[two heads]{d} & Q \arrow[equal]{d} \\
Y \arrow[tail]{r} & P \arrow[two heads]{r} & Q
\end{tikzcd}$$

\

\noindent \underline{Claim 3}. Next suppose $W  \rightarrowtail W' \twoheadrightarrow  Z$ is exact with $W, W' \in \class{W}$. We need to prove that $Z \in \class{W}$ too. Start by writing $W  \rightarrowtail \tilclass{R} \twoheadrightarrow  \tilclass{Q}$ and again forming a pushout diagram as shown.
$$\begin{tikzcd}
W \arrow[tail]{r} \arrow[tail]{d} & W' \arrow[two heads]{r} \arrow[tail]{d} & Z \arrow[equal]{d} \\
\tilclass{R} \arrow[tail]{r} \arrow[two heads]{d} & P \arrow[two heads]{r} \arrow[two heads]{d} & Z  \\
\tilclass{Q} \arrow[equal]{r} & \tilclass{Q}
\end{tikzcd}$$ Since we have shown $\class{W}$ is closed under extensions we get $P \in \class{W}$. So now we can write $P  \rightarrowtail \tilclass{R}' \twoheadrightarrow  \tilclass{Q}'$. We now take yet another pushout, the pushout of $\tilclass{R}' \leftarrowtail P \twoheadrightarrow Z$ to get yet another diagram as shown.
$$\begin{tikzcd}
\tilclass{R} \arrow[equal]{r} \arrow[tail]{d} & \tilclass{R} \arrow[tail]{d} &  \\
P \arrow[tail]{r} \arrow[two heads]{d} & \tilclass{R}' \arrow[two heads]{r} \arrow[two heads]{d} & \tilclass{Q}' \arrow[equal]{d} \\
Z \arrow[tail]{r} & L \arrow[two heads]{r} & \tilclass{Q}'
\end{tikzcd}$$
(By~\cite[Proposition~2.12]{buhler-exact categories} we get that the lower left square is bicartesian.) Now since $(\class{Q}, \tilclass{R})$ is an hereditary cotorsion pair we get $L \in \tilclass{R}$. So the short exact sequence $Z \rightarrowtail L \twoheadrightarrow \tilclass{Q}'$ finishes the proof of the claim.

A dual argument shows that if $X  \rightarrowtail W' \twoheadrightarrow  W$ is exact with $W, W' \in \class{W}$ then so is $X$. (Use the \emph{other} characterization of $\class{W}$ and pullbacks.) This completes the proof that $\class{W}$ is thick.

\

\noindent ($(\class{Q},\class{W},\class{R})$ is a Hovey triple.) We just need to show $\class{Q} \cap \class{W} = \tilclass{Q}$ and $\class{W} \cap \class{R} = \tilclass{R}$. The proof of each is similar, so we will just show $\class{W} \cap \class{R} = \tilclass{R}$. Here, we clearly have $\class{W} \cap \class{R} \supseteq \tilclass{R}$, so we just need to show $\class{W} \cap \class{R} \subseteq \tilclass{R}$. So let $X \in \class{W} \cap \class{R}$. Using that $X$ is in $\class{W}$ we may write a short exact sequence $X  \rightarrowtail \widetilde{R} \twoheadrightarrow  \widetilde{Q}$. But since $X \in \class{R}$ and $(\tilclass{Q},\class{R})$ is a cotorsion pair, this sequence must split. So $X$ is a retract of $\widetilde{R}$ and so must be in $\tilclass{R}$. This completes the proof that $(\class{Q},\class{W},\class{R})$ is a Hovey triple. The uniqueness of $\class{W}$ follows from a general fact: The class of trivial objects in a Hovey triple is always unique by~\cite[Proposition~3.2]{gillespie-recollements}.
\end{proof}

\end{document}